%% file: m3-II-5.tex
\input m3-macs

\pageno=255

\tinfo II.5.255-262

\SetTFLinebox{\gtp }
\SetFLinebox{\gtv3 }
\SetHLinebox{\issn}

\H 5. Harmonic analysis on algebraic groups over\\
two-dimensional local fields of equal characteristic

Mikhail Kapranov

\SetAuthorHead{M. Kapranov}
\SetTitleHead{Part II. Section 5. Harmonic analysis on algebraic groups \qquad\qquad}

In this section we review the main parts of a recent  work \cite{4}
on harmonic analysis on algebraic groups over 
two-dimensional local fields.

\HH 5.1. Groups and buildings

Let $K$ ($K=K_2$ whose residue field is $K_1$ whose residue field is $K_0$, see the notation in section~1 of Part~I) 
be a two-dimensional local field of equal characteristic.
Thus $K_2$  is isomorphic to the Laurent
series field $K_1((t_2))$ over $K_1$. It is convenient to think of elements
of $K_2$ as (formal) loops over $K_1$. Even in the case
where $\chr(K_1)=0$, it is still convenient
to think of elements of $K_1$ as (generalized) loops over $K_0$
so that $K_2$ consists of double loops. 

Denote the residue map $\Cal O_{K_2}\to K_1$ by $p_2$
and  the residue map $\Cal O_{K_1}\to K_0$ by $p_1$.
Then the ring of integers $O_K$ of $K$ as
of a two-dimensional local field (see subsection 1.1 of Part I) coincides with
$p_2^{-1}(\Cal O_{K_1})$.

Let $G$ be a split simple simply connected algebraic group over $\Bbb Z$
(e.g. $G=SL_2$).
Let $T\subset B\subset G$ be a fixed maximal torus and Borel subgroup of $G$; 
put $N=[B,B]$, and let $W$ be the Weyl group of $G$.
All of them are viewed as group schemes.

Let $L=\Hom(\Bbb G_m,T)$ be the coweight  lattice of $G$;
the Weyl group acts on $L$.

Recall that $I(K_1)=p_1^{-1}(B(\Bbb F_q))$ is called an Iwahori
subgroup of $G(K_1)$ and $T(\Cal O_{K_1})N(K_1)$ can be
seen as the ``connected component of unity'' in
$B(K_1)$. The latter name is explained naturally if we think
of elements of $B(K_1)$ as being loops with values in $B$.

\df Definition

Put
$$\aligned
&D_0=p_2^{-1}p_1^{-1}(B(\Bbb F_q))\i G(O_K), \\
&D_1=p_2^{-1}(T(\Cal O_{K_1})N(K_1))\i G(O_K), \\
&D_2=T(\Cal O_{K_2})N(K_2)\i G(K).
\endaligned $$
\enddf

Then $D_2$ can be seen as the ``connected component of unity''
 of $B(K)$ when $K$ is viewed as a two-dimensional local field,
$D_1$ is a (similarly understood)
connected component of an Iwahori subgroup of $G(K_2)$, and 
$D_0$ is called a double Iwahori subgroup of $G(K)$.

A choice of a system of local parameters $t_1,t_2$ of $K$ determines
the identification of the group $K^*/O_K^*$ with $\Bbb Z\oplus \Bbb Z$  
and identification $L\oplus L$ with  $L\otimes (K^*/O_K^*)$.

We have 
 an embedding of $L\otimes (K^*/O_K^*)$ into $T(K)$
which takes $a\otimes (t_1^j t_2^j)$, $i, j\in\Bbb Z$,
to the value on $t_1^i t_2^j$ of the 1-parameter
subgroup in $T$ corresponding to $a$.

Define the action of $W$ on $L\otimes (K^*/O_K^*)$ as the product
of the standard action on $L$ and the  trivial action on
 $K^*/O_K^*$.
The semidirect product
$$\wh{\wh{W}}=(L\otimes K^*/O_K^*)\sdp W$$
is called the {\it double affine Weyl group} of $G$.

A (set-theoretical) lifting  of $W$ into $G(O_K)$ determines
a lifting of  $\wh{\wh{W}}$ into $G(K)$.

\th Proposition

For every $i,j\in\{0,1,2\}$ there is a disjoint decomposition
$$G(K)=\bigcupp_{w\in \wh{\wh{W}}} D_iwD_j.$$
The identification $D_i\backslash G(K)/D_j$ with $\wh{\wh{W}}$
doesn't depend on the choice of liftings.
\endth

\pf Proof

Iterated application of the Bruhat, Bruhat--Tits and 
Iwasawa decompositions to the local fields
$K_2$, $K_1$.
\endpf

For the Iwahori subgroup $I(K_2)=p_2^{-1}(B(K_1))$ of $G(K_2)$  
the homogeneous space $G(K)/I(K_2)$ is the ``affine flag variety''
of $G$, see \cite{5}. It has a canonical structure of an ind-scheme,
in fact, it is an inductive limit of projective algebraic
varieties over $K_1$ (the closures of the affine Schubert cells). 

Let $B(G,K_2/K_1)$ be the {\it Bruhat--Tits building} associated to $G$ and 
the field $K_2$. Then the space $G(K)/I(K_2)$ is a $G(K)$-orbit on the set of 
flags of type (vertex, maximal cell) in the building.
For every vertex $v$ of $B(G,K_2/K_1)$
its locally finite Bruhat--Tits building
$\beta_v$ isomorphic to $B(G,K_1/K_0)$
can be viewed as a ``microbuilding'' of the {\it double Bruhat--Tits
building} $B(G,K_2/K_1/K_0)$ of $K$ as a two-dimensional local field
constructed by Parshin (\cite{7}, see also section~3 of Part~II).
Then the set $G(K)/D_1$ is identified naturally
with the set of all the horocycles
$\{w\in \beta_v: d(z,w)=r\}$, $z\in \partial\beta_v$
of the microbuildings $\beta_v$ 
(where the ``distance'' $d(z,\,\,)$ is viewed as an element of
a natural $L$-torsor).
The fibres of the projection $G(K)/D_1\to G(K)/I(K_2)$
are $L$-torsors.

\vskip .3cm

\HH 5.2. The central extension and the affine Heisenberg--Weyl group

According to  the work of Steinberg, Moore and Matsumoto \cite{6} 
developed by Brylinski and Deligne \cite{1} 
there is a central extension
$$1 \to K_1^*\to \Gamma \to G(K_2)\to 1$$
associated to
the tame symbol $K_2^*\times K_2^* \to K_1^*$
for the couple $(K_2,K_1)$ (see subsection 6.4.2 of Part~I for the general 
definition of the tame symbol).

\th Proposition

This extension splits over every $D_i$, $0\le i\le 2$.
\endth

\pf Proof

Use Matsumoto's explicit construction of the central extension.
\endpf

Thus, there are identifications of every $D_i$ with a subgroup of $\Gamma$.
Put
$$\Delta_i=\Cal O_{K_1}^*D_i\i \Gamma, \qquad \Xi=\Gamma/\Delta_1.$$

The minimal integer scalar product $\Psi$ on $L$ 
and the composite of the tame symbol 

\noindent $K_2^*\times K_2^*\to K_1^*$ and the discrete valuation
$v_{K_1}\colon K^*\to \Bbb Z$
induces a $W$-invariant skew-symmetric pairing $L\otimes K^*/O_K^*
 \times
L\otimes K^*/O_K^*\to \Bbb Z$.
Let
$$1\to \Bbb Z\to \Cal L\to  L\otimes K^*/O_K^*\to 1$$
be the  central extension
whose commutator pairing corresponds 
to  the latter skew-symmetric pairing. The group
$\Cal L$ is called the {\it Heisenberg group}.

\df Definition

The semidirect product
$$\wt{W}=\Cal L\sdp W$$
is called the {\it double affine Heisenberg--Weyl} group of $G$.
\enddf

\th Theorem

The group $\wt{W}$ is isomorphic to $L_{\aff}\sdp \wh{W}$
where $L_{\aff}=\Bbb Z\oplus L$, 
$\wh{W}=L\sdp W$ and
$$w\circ (a,l')=(a,w(l)),\quad l\circ (a,l')=(a+\Psi(l,l'),l'),\qquad
w\in W,\quad l,l'\in L,\quad a\in \Bbb Z.$$
 For every $i,j\in\{0,1,2\}$ there is a disjoint union
$$\Gamma=\bigcupp_{w\in\wt{W}} \Delta_iw\Delta_j$$
and the identification $\Delta_i\backslash \Gamma/\Delta_j$ with $\wt{W}$
is canonical.
\endth

\HH 5.3. Hecke algebras in the classical setting

Recall that for a locally compact group $\Gamma$ and its compact subgroup 
$\Delta$
the Hecke algebra $\Cal H(\Gamma,\Delta)$ can be defined as the algebra 
 of compactly supported  double $\Delta$-invariant continuous
functions of $\Gamma$ with the operation given by the convolution
with respect to the Haar measure on $\Gamma$.
For  $C=\Delta\gamma \Delta\in \Delta\backslash \Gamma/\Delta$
the Hecke correspondence $\Sigma_C=\{(\alpha \Delta, \beta \Delta): 
\alpha\beta^{-1}\in C\}$ is a $\Gamma$-orbit of $(\Gamma/\Delta)\times(\Gamma/
\Delta)$.

For $x\in \Gamma/\Delta$ put $\Sigma_C(x)=
\Sigma_C\cap (\Gamma/\Delta)\times
\{x\}$. Denote the projections of $\Sigma_C$
 to the first and second component
by $\pi_1$ and $\pi_2$.

Let $\Cal F(\Gamma/\Delta)$ be the space of continuous functions
$\Gamma/\Delta\to\Bbb C$.
The  operator 
$$\tau_C\colon \Cal F(\Gamma/\Delta)\to  \Cal F(\Gamma/\Delta),
\quad f\to {\pi_2}_*\pi_1^*(f)
$$ is called the {\it Hecke operator} associated to $C$. 
Explicitly, 
$$(\tau_Cf)(x)=\int_{y\in\Sigma_C(x)} f(y)d\mu_{C,x},$$
where $\mu_{C,x}$ is the $\Stab(x)$-invariant 
measure induced by the Haar measure.
Elements of the Hecke algebra $\Cal H(\Gamma,\Delta)$
can be viewed as ``continuous''  
linear combinations of the operators $\tau_C$, i.e., integrals
of the form $\int \phi(C) \tau_C dC$ where $dC$ is some measure
on $\Delta\backslash\Gamma/\Delta$
and $\phi$ is a continuous function with compact support.
 If the group $\Delta$ is also
open (as is usually the case in the $p$-adic situation), then
$\Delta\backslash\Gamma/\Delta$ is discrete and $\Cal H(\Gamma, \Delta)$
consists of finite linear combinations of the $\tau_C$.

\HH 5.4. The regularized  Hecke algebra $\Cal H(\Gamma,\Delta_1)$

Since the two-dimensional local field
$K$ and the ring $O_K$ are not locally compact, the  approach
of the previous subsection would work only after a new appropriate integration theory is available.

The aim of this subsection is to make sense of the Hecke algebra
$\Cal H(\Gamma,\Delta_1)$.

Note that the fibres of the projection $\Xi=\Gamma/\Delta_1\to G(K)/I(K_2)$
are $ L_{\aff}$-torsors and
$G(K)/I(K_2)$ is the inductive limit of compact (profinite) spaces,
so $\Xi$ can be considered as an object of the category
$\Cal F_1$ defined in subsection~1.2 of the paper of Kato in this volume.

Using Theorem of 5.2 for $i=j=1$ we introduce:

\df Definition

For $(w,l)\in \wt{W}=L_{\aff}\sdp \wh{W}$
denote by $\Sigma_{w,l}$ the Hecke correspondence 
(i.e., the $\Gamma$-orbit
of $\Xi\times \Xi$) associated to $(w,l)$.
For $\xi\in \Xi$ put
$$\Sigma_{w,l}(\xi)=\{\xi':(\xi,\xi')\in \Sigma_{w,l}\}.$$
\enddf

The stabilizer $\Stab(\xi)\le \Gamma$ acts
transitively on $\Sigma_{w,l}(\xi)$.

\th Proposition

$\Sigma_{w,l}(\xi)$ is an affine space over $K_1$ of dimension
equal to the length of $w\in \wh{W}$.
The space of compex valued Borel measures on $\Sigma_{w,l}(\xi)$
is 1-dimensional.
A choice of a $\Stab(\xi)$-invariant measure $\mu_{w,l,\xi}$
on $\Sigma_{w,l}(\xi)$
determines a measure $\mu_{w,l,\xi'}$ on $\Sigma_{w,l}(\xi')$ for every $\xi'$.
\endth

\df Definition

For a continuous function $f\colon \Xi\to \Bbb C$
put
$$(\tau_{w,l}f)(\xi)=\int_{\eta\in \Sigma_{w,l}(\xi)}
 f(\eta)d\mu_{w,l,\xi}.$$

\enddf

Since the domain of the integration is not compact,
the integral may diverge. As a first step, we define
the space of functions on which the integral makes
sense. Note that $\Xi$ can be regarded as an $L_{\aff}$-torsor
over the ind-object $G(K)/I(K_2)$ in the category $\pro (C_0)$,
i.e., a compatible system of $L_{\aff}$-torsors $\Xi_\nu$ over the
affine Schubert varieties $Z_\nu$ forming an exhaustion of $G(K)/I(K_1)$.
Each $\Xi_\nu$ is a locally compact space and $Z_\nu$
is a compact space. In particular,  we can form the
space $\Cal F_0 (\Xi_\nu)$ of locally constant
complex  valued functions on $\Xi_\nu$
whose support is compact (or, what is the same, proper
with respect to the projection to $Z_\nu$).
Let $\Cal F(\Xi_\nu)$ be the space of all
locally constant complex functions on $\Xi_\nu$. Then we define
$\Cal F_0(\Xi) =$"$\varprojlim$"$ \Cal F_0(\Xi_\nu)$
and $\Cal F(\Xi) =  $"$\varprojlim $"$ \Cal F(\Xi_\nu)$.
 They are
 pro-objects in the category of  vector spaces. 
In fact, because of the action of $L_{\aff}$ and its
group algebra $\Bbb C[L_{\aff}]$ on $\Xi$,  the spaces
$\Cal F_0(\Xi), \Cal F(\Xi)$ are naturally pro-objects
in the category of $\Bbb C[L_{\aff}]$-modules.

\th Proposition

If $f=(f_\nu)\in \Cal F_0(X)$
then $\text{\rm Supp}(f_\nu)\cap \Sigma_{w,l}(\xi)$
is compact for every $w,l,\xi, \nu$ and the integral above
converges.
Thus, there is a well defined  Hecke operator
$$\tau_{w,l}\colon \Cal F_0(\Xi)\to\Cal F(\Xi)$$
which is an element of
$\text{\rm Mor(pro(Mod${}_{\Bbb C\,[L_{\aff}]}$))}$. 
In particular, $\tau_{w,l}$ is the shift by $l$
and $\tau_{w,l+l'}=\tau_{w,l'}\tau_{e,l}$.
\endth

Thus we get Hecke operators as operators from one
(pro-)vector space to another, bigger one. This does not
yet allow to compose the $\tau_{w, l}$. Our next step
is to consider certain infinite linear combinations
of the $\tau_{w, l}$. 

Let $T^\vee_{\aff} = \text{Spec}(\Bbb C[L_{\aff}])$
be the ``dual affine torus'' of $G$. A function
with finite support on $L_{\aff}$ can be viewed
as the collection of coefficients of a polynomial, i.e., of
 an element of $\Bbb C[L_{\aff}]$  as a regular function on
$T^\vee_{\aff}$.  Further, let $Q\i L_{\aff}\otimes\Bbb R$
be a strictly convex cone with apex $0$. A function
on $L_{\aff}$ with support in $Q$ can be viewed as
the collection of coefficients of a formal power
series, and such series form a ring containing
$\Bbb C[L_{\aff}]$. On the level of functions
the ring operation is the convolution. 
 Let $\Cal F_Q(L_{\aff})$ be the
space of functions whose support is contained
in some translation of $Q$. It is a ring with
respect to convolution.

 Let $\Bbb C(L_{\aff})$ be the field
of rational functions on $T^\vee_{\aff}$.
 Denote by 
$F_Q^{\text{rat}}(L_{\aff})$ the subspace in $F_Q(L_{\aff})$
consisting of functions whose corresponding
formal power series are expansions of rational functions
on $T^\vee_{\aff}$.

If $A$ is any $L_{\aff}$-torsor (over a point), then $\Cal F_0(A)$
is an (invertible) module over $\Cal F_0(L_{\aff})=\Bbb C[L_{\aff}]$
and we can define the spaces $\Cal F_Q(A)$ and
$\Cal F_Q^{\text{rat}}(A)$ which will be modules
over the corresponding rings for $L_{\aff}$.
We also write $\Cal F^{\text{rat}}(A) =
 \Cal F_0(A)\otimes_{\Bbb C[L_{\aff}]}
\Bbb C(L_{\aff})$. 

We then extend the above concepts ``fiberwise'' to torsors
over compact spaces (objects of $\pro(C_0)$)
and to torsors over objects of $\ind(\pro(C_0))$
such as $\Xi$. 

Let $w\in \widehat W$. We denote by $Q(w)$ the image under $w$
of the cone of dominant affine coweights in $L_{\aff}$.

\th Theorem

The action of the Hecke operator $\tau_{w,l}$ takes $\Cal F_0(\Xi)$
into $\Cal F_{Q(w)}^{\text{rat}}(\Xi)$. These operators
extend to operators 
$$\tau_{w,l}^{\text{rat}}: \Cal F^{\text{rat}}(\Xi)\to 
\Cal F^{\text{rat}}(\Xi).$$

\endth

Note that the action of $\tau_{w,l}^{\text{rat}}$ involves
a kind of regularization procedure, which is hidden
in the identification of the $\Cal F_{Q(w)}^{\text{rat}}(\Xi)$
for different $w$,  with subspaces of the same space
$\Cal F^{\text{rat}}(\Xi)$. In practical terms, this
involves summation of a series to a rational function
and re-expansion in a different domain.

\smallskip

Let $\Cal H_{\text{pre}}$ be the space of finite linear combinations
  $\sum_{w,l} a_{w,l}\tau_{w,l}$. This is not yet an algebra,
but only a $\Bbb C[L_{\aff}]$-module. 
Note that elements of $\Cal H_{\text{pre}}$ can be written as
 finite linear combinations $\sum_{w\in\widehat W} f_w(t) \tau_w$
where $f_w(t) = \sum_l a_{w,l} t^l$, $t\in T^\vee_{\aff}$, is the
polynomial in $\Bbb C[L_{\aff}]$ corresponding to the collection of the 
$a_{w,l}$.
This makes the $\Bbb C[L_{\aff}]$-module structure clear. 
Consider the tensor product
$$\Cal H_{\text{rat}} = \Cal H_{\text{pre}}
\otimes_{\Bbb C\,[L_{\aff}]}\Bbb C(L_{\aff}).$$
Elements of this space can be considered as finite linear
combinations  $\sum_{w\in\widehat W} f_w(t) \tau_w$
where $f_w(t)$ are now rational functions. By expanding
rational functions in power series, we can
consider the above elements as certain infinite linear combinations
of the $\tau_{w,l}$.

\th Theorem

The space $\Cal H_{\text{rat}}$ has a natural algebra structure
and this algebra acts in the space  $\Cal F^{\text{rat}}(\Xi)$,
 extending the action of the $\tau_{w,l}$ defined above.
\endth

The operators associated to $\Cal H_{\text{rat}}$ can be
viewed as certain integro-difference operators,
because their action involves integration (as in the
definition of the $\tau_{w,l}$) as well as inverses
of linear combinations of shifts by elements
of $L$ (these combinations act as difference operators).

 \df Definition

The regularized Hecke algebra $\Cal H(\Gamma, \Delta_1)$
is, by definition, the subalgebra in $\Cal H_{\text{rat}}$
consisting of elements whose action in $\Cal F_{\text{rat}}(\Xi)$
preserves the subspace $\Cal F_0(\Xi)$. 

\enddf

\HH 5.5. The Hecke algebra and  the Cherednik algebra

In \cite{2} I. Cherednik introduced the so-called
  double affine Hecke algebra $\text{Cher}_q$
 associated
to the root system of $G$.
As shown by V. Ginzburg, E. Vasserot and the author \cite{3},
$\text{Cher}_q$ can be thought as consisting of finite
linear combinations $\sum_{w\in\widehat W_{\text{ad}}} f_w(t) [w]$
where $W_{\text{ad}}$ is the affine Weyl group
of the adjoint quotient $G_{\text{ad}}$ of $G$ (it contains $\widehat W$)
and $f_w(t)$ are rational functions on $T^{\vee}_{\aff}$
satisfying certain residue conditions. We define
the modified Cherednik algebra $\ddot{H}_q$
to be the subalgebra in $\text{Cher}_q$ consisting of
linear combinations as above, but going over $\widehat W\i \widehat W_{
\text{ad}}$.

\th Theorem

The regularized Hecke algebra $\Cal H(\Gamma, \Delta_1)$ is isomorphic
to the modified Cherednik algebra $\ddot{H}_q$.
In particular, there is a natural action of $\ddot{H}_q$
on $\Cal F_0(\Xi)$ by integro-difference operators.
\endth
\pf Proof

Use the principal series intertwiners and a version of Mellin transform.
The information on the poles of the intertwiners matches
exactly the residue conditions introduced in \cite{3}. 
\endpf

\rk Remark

 The only reason we needed to assume that
the 2-dimensional local field $K$ has equal characteristic
was because we used the fact that the quotient
$G(K)/I(K_2)$ has a structure of an inductive limit of
projective algebraic varieties over $K_1$. In fact, we really
use only a weaker structure: that of an inductive
limit of  profinite topological spaces
 (which are, in this case, the sets of $K_1$-points of
 affine Schubert varieties over $K_1$).
This structure is available for any 2-dimensional local field,
although there seems to be no reference for it in the literature.
Once this foundational matter is established, all the constructions will
go through for any 2-dimensional local field.

\endrk 

\vskip 1cm

\Bib References

\rf{1} J.-L. Brylinski and P. Deligne, Central extensions of  reductive
groups by $\Cal K_2$, preprint of IAS, Princeton, available
from P. Deligne's home page at
 $<$www.math.ias.edu$>$.

\rf{2} I. Cherednik, Double affine Hecke algebras and Macdonald's
conjectures, {Ann. Math.} {141} (1995), 191--216. 

\rf{3} V. Ginzburg and M. Kapranov and E. Vasserot, Residue construction
of Hecke algebras, {Adv. in Math.} {128} (1997), 1--19.

\rf{4} M. Kapranov, Double affine Hecke algebras and 2-dimensional
local fields, preprint math.AG/9812021, to appear in 
Journal of the AMS. 

\rf{5} G. Lusztig,
  Singularities,
character
formula and $q$-analog of weight multiplicity,
 {Asterisque} {
101-102} (1983), 208--222.

\rf{6} H. Matsumoto, Sur les sous-groupes arithm\'etiques
des groupes semi-simples d\'eploy\'es, {Ann. ENS}, {2} (1969),
1--62. 

\rf{7} A. N. Parshin,  
Vector bundles and arithmetic groups I: The higher Bruhat-Tits tree,
 {Proc. Steklov Inst. Math.} {208} (1995), 212--233,
preprint alg-geom/9605001.

\endBib 

\Coordinates

Department of Mathematics \  University of Toronto

Toronto  M5S 3G3 \ Canada 

E-mail: kapranov\@math.toronto.edu

\endCoordinates

\vfill
\pagebreak

\end

%% file: m3-macs.tex
\expandafter\ifx\csname mthreemacsloaded\endcsname\relax\else \fi

\magnification1100
\input amstex


 \catcode`\@=11
 \let\wlog@ld\wlog
 \def\wlog#1{\relax}

 \newif\ifIN@
 \def\m@rker{\m@@rker}
 \def\IN@{\expandafter\INN@\expandafter}
 \long\def\INN@0#1@#2@{\long\def\NI@##1#1##2##3\ENDNI@
    {\ifx\m@rker##2\IN@false\else\IN@true\fi}%
     \expandafter\NI@#2@@#1\m@rker\ENDNI@}
  \newtoks\Initialtoks@  \newtoks\Terminaltoks@
  \def\SPLIT@{\expandafter\SPLITT@\expandafter}
  \def\SPLITT@0#1@#2@{\def\TTILPS@##1#1##2@{%
     \Initialtoks@{##1}\Terminaltoks@{##2}}\expandafter\TTILPS@#2@}
  \newtoks\Trimtoks@

 \def\ForeTrim@{\expandafter\ForeTrim@@\expandafter}
 \def\ForePrim@0 #1@{\Trimtoks@{#1}}
 \def\ForeTrim@@0#1@{\IN@0\m@rker. @\m@rker.#1@%
     \ifIN@\ForePrim@0#1@%
     \else\Trimtoks@\expandafter{#1}\fi}
 
  \def\Trim@0#1@{%
      \ForeTrim@0#1@%
      \IN@0 @\the\Trimtoks@ @%
        \ifIN@
             \SPLIT@0 @\the\Trimtoks@ @\Trimtoks@\Initialtoks@
             \IN@0\the\Terminaltoks@ @ @%
                 \ifIN@
                 \else \Trimtoks@ {FigNameWithSpace}%
                 \fi
        \fi
      }

  \font\titlebold=cmbx12 scaled 1200
  \font\twelvebold=cmbx12
  \font\tenbold=cmbx10
  \font\ninebold=cmbx9
  \font\sevenbold=cmbx7
  \font\fivebold=cmbx5

  \input amssym.def \input amssym
     \font\titlemsa=msam10 at 14.4pt
     \font\titlemsb=msbm10 at 14.4pt
     \font\titleeufm=eufm10 at 14.4pt
     \font\twelvemsa=msam10 scaled 1200
     \font\twelvemsb=msbm10 scaled 1200
     \font\twelveeufm=eufm10 scaled 1200
     \font\ninemsa=msam9
     \font\ninemsb=msbm9
     \font\nineeufm=eufm9

   \ifx\cyrfam\undefined
   \else
     \immediate\write16{}%
     \message{ !!! cyr fonts already defined. !!! }
     \message{ --- edit out superfluous font defs? }
   \fi
   \newfam\cyrfam
       \font\titlecyr=wncyr10 scaled 1440 
       \font\twelvecyr=wncyr10 scaled 1200
       \font\tencyr=wncyr10
       \font\ninecyr=wncyr9
       \font\sevencyr=wncyr7
       \font\sixcyr=wncyr6

   \newfam\eusmfam
       \font\titleeusm=eusm10 scaled 1440
       \font\twelveeusm=eusm10 scaled 1200
       \font\teneusm=eusm10
       \font\nineeusm=eusm9
       \font\seveneusm=eusm7
       
       \font\fiveeusm=eusm5

\let\Cal\cal

    \font\ninemrm=cmr9 
    \font\ninei=cmmi9
    \font\ninesy=cmsy9 
    \skewchar\ninei='177
    \skewchar\ninesy='60

  \font\twelvemrm=cmr10 at 12pt 
  \font\twelvei=cmmi10 at 12pt
  \font\twelvesy=cmsy10 at 12pt

  \font\titlemrm=cmr10 at 14.4pt 
  \font\titlei=cmmi10 at 14.4pt
  \font\titlesy=cmsy10 at 14.4pt


  \def\Smallfonts{\ninepoint}

  \def\Hfont{\titlepoint\bf}
  \def\Authorfont{\twelvepoint\it}
  \def\HHfont{\twelvepoint\bf}
  \def\HHHfont{\bf}
  \def\Bibfont{\tenbf}
  \def\Coordfont{\nineit }

  \def \thfont {\bf }
  \def \pffont {\it\itSpacing }
  \def \rkfont {\bf }
  \def \dffont {\bf }
  \def \egfont {\bf }

 \def\ninepoint{%
  \def\rm{\fam0\ninerm}%
    \textfont0=\ninemrm  \scriptfont0=\sevenrm  \scriptscriptfont0=\fiverm
    \textfont1=\ninei    \scriptfont1=\seveni   \scriptscriptfont1=\fivei
  \def\mit{\fam1\ninei}%
  \def\oldstyle{\fam1\ninei}%
    \textfont2=\ninesy   \scriptfont2=\sevensy  \scriptscriptfont2=\fivesy
    \textfont3=\tenex    \scriptfont3=\tenex    \scriptscriptfont3=\tenex
  \def\it{\fam\itfam\nineit}%
    \textfont\itfam=\nineit
  \def\bf{\ifmmode\fam\bffam\else\ninebf\fi}%
    \textfont\bffam=\ninebold 
    \scriptfont\bffam=\sevenbold 
    \scriptscriptfont\bffam=\fivebold%
  \def\msa{\fam\msafam\ninemsa}%
    \textfont\msafam=\ninemsa 
    \scriptfont\msafam=\sevenmsa
    \scriptscriptfont\msafam=\fivemsa%
  \def\msb{\fam\msbfam\ninemsb}%
    \textfont\msbfam=\ninemsb%
    \scriptfont\msbfam=\sevenmsb%
    \scriptscriptfont\msbfam=\fivemsb%
  \def\eufm{\fam\eufmfam\nineeufm}%
    \textfont\eufmfam=\nineeufm
    \scriptfont\eufmfam=\seveneufm
    \scriptscriptfont\eufmfam=\fiveeufm
   \def\eusm{\fam\eusmfam\nineeusm}%
     \textfont\eusmfam=\nineeusm
     \scriptfont\eusmfam=\seveneusm
     \scriptscriptfont\eusmfam=\fiveeusm
   \def\cyr{\fam\cyrfam\ninecyr}%
     \textfont\cyrfam=\ninecyr
     \scriptfont\cyrfam=\sevencyr
     \scriptscriptfont\cyrfam=\sixcyr
  \setbox\strutbox=\hbox{\vrule
      height7pt depth3pt width0pt}%
   \baselineskip=10.8pt\rm}

 \let\eightpoint\ninepoint 

 \def\tenpoint{%
  \def\rm{\fam0\tenrm}%
    \textfont0=\tenmrm \scriptfont0=\sevenrm \scriptscriptfont0=\fiverm%
  \def\mit{\fam1\teni}%
  \def\oldstyle{\fam1\teni}%
    \textfont1=\teni   \scriptfont1=\seveni  \scriptscriptfont1=\fivei%
    \textfont2=\tensy  \scriptfont2=\sevensy \scriptscriptfont2=\fivesy%
    \textfont3=\tenex  \scriptfont3=\tenex   \scriptscriptfont3=\tenex%
  \def\it{\fam\itfam\tenit}%
    \textfont\itfam=\tenit%
  \def\bf{\ifmmode\fam\bffam\else\tenbf\fi}%
    \textfont\bffam=\tenbold
    \scriptfont\bffam=\sevenbold%
    \scriptscriptfont\bffam=\fivebold%
  \def\msa{\fam\msafam\tenmsa}%
    \textfont\msafam=\tenmsa%
    \scriptfont\msafam=\sevenmsa%
    \scriptscriptfont\msafam=\fivemsa%
  \def\msb{\fam\msbfam\tenmsb}%
    \textfont\msbfam=\tenmsb%
    \scriptfont\msbfam=\sevenmsb%
    \scriptscriptfont\msbfam=\fivemsb%
  \def\eufm{\fam\eufmfam\teneufm}%
   \textfont\eufmfam=\teneufm
   \scriptfont\eufmfam=\seveneufm
   \scriptscriptfont\eufmfam=\fiveeufm
   \def\eusm{\fam\eusmfam\teneusm}%
    \textfont\eusmfam=\teneusm
    \scriptfont\eusmfam=\seveneusm
    \scriptscriptfont\eusmfam=\fiveeusm
   \def\cyr{\fam\cyrfam\tencyr}%
    \textfont\cyrfam=\tencyr
    \scriptfont\cyrfam=\sevencyr
    \scriptscriptfont\cyrfam=\sixcyr
  \setbox\strutbox=\hbox{\vrule %
      height8.5pt depth3.5ptwidth0pt}%
  \baselineskip=\StdBaselineskip\rm}

 \def\twelvepoint{%
  \def\rm{\fam0\twelverm}%
    \textfont0=\twelvemrm \scriptfont0=\tenmrm \scriptscriptfont0=\sevenrm
    \textfont1=\twelvei   \scriptfont1=\teni   \scriptscriptfont1=\seveni
  \def\mit{\fam1\twelvei}%
  \def\oldstyle{\fam1\twelvei}%
    \textfont2=\twelvesy  \scriptfont2=\tensy  \scriptscriptfont2=\sevensy
    \textfont3=\tenex  \scriptfont3=\tenex  \scriptscriptfont3=\tenex
  \def\it{\fam\itfam\twelveit}%
    \textfont\itfam=\twelveit
  \def\bf{\ifmmode\fam\bffam\else\twelvebf\fi}%
    \textfont\bffam=\twelvebold
    \scriptfont\bffam=\tenbold%
    \scriptscriptfont\bffam=\sevenbold%
  \def\msa{\fam\msafam\twelvemsa}%
    \textfont\msafam=\twelvemsa%
    \scriptfont\msafam=\tenmsa%
    \scriptscriptfont\msafam=\sevenmsa%
  \def\msb{\fam\msbfam\twelvemsb}%
    \textfont\msbfam=\twelvemsb%
    \scriptfont\msbfam=\tenmsb%
    \scriptscriptfont\msbfam=\sevenmsb%
  \def\eufm{\fam\eufmfam\twelveeufm}%
   \textfont\eufmfam=\twelveeufm
   \scriptfont\eufmfam=\teneufm
   \scriptscriptfont\eufmfam=\seveneufm
   \def\eusm{\fam\eusmfam\twelveeusm}%
    \textfont\eusmfam=\twelveeusm
    \scriptfont\eusmfam=\teneusm
    \scriptscriptfont\eusmfam=\seveneusm
   \def\cyr{\fam\cyrfam\tencyr}%
    \textfont\cyrfam=\twelvecyr
    \scriptfont\cyrfam=\tencyr
    \scriptscriptfont\cyrfam=\sevencyr
  \setbox\strutbox=\hbox{\vrule
      height10.2pt depth4.55pt width0pt}%
  \baselineskip=14pt\rm}

 \def\titlepoint{%
    \textfont0=\titlemrm \scriptfont0=\twelvemrm \scriptscriptfont0=\tenmrm
    \textfont1=\titlei   \scriptfont1=\twelvei   \scriptscriptfont1=\teni
  \def\mit{\fam1\titlei}%
  \def\oldstyle{\fam1\titlei}%
    \textfont2=\titlesy  \scriptfont2=\twelvesy  \scriptscriptfont2=\tensy
    \textfont3=\tenex
    \scriptfont3=\tenex
    \scriptscriptfont3=\tenex
  \def\it{\fam\itfam\titleit}%
    \textfont\itfam=\titleit
  \def\bf{\ifmmode\fam\bffam\else\titlebf\fi}%
    \textfont\bffam=\titlebold
    \scriptfont\bffam=\twelvebold%
    \scriptscriptfont\bffam=\tenbold%
  \def\msa{\fam\msafam\titlemsa}%
    \textfont\msafam=\titlemsa%
    \scriptfont\msafam=\twelvemsa%
    \scriptscriptfont\msafam=\tenmsa%
  \def\msb{\fam\msbfam\titlemsb}%
    \textfont\msbfam=\titlemsb%
    \scriptfont\msbfam=\twelvemsb%
    \scriptscriptfont\msbfam=\tenmsb%
  \def\eufm{\fam\eufmfam\titleeufm}%
    \textfont\eufmfam=\titleeufm
    \scriptfont\eufmfam=\twelveeufm
    \scriptscriptfont\eufmfam=\teneufm
   \def\eusm{\fam\eusmfam\titleeusm}%
     \textfont\eusmfam=\titleeusm
     \scriptfont\eusmfam=\twelveeusm
     \scriptscriptfont\eusmfam=\teneusm
   \def\cyr{\fam\cyrfam\tencyr}%
    \textfont\cyrfam=\titlecyr
    \scriptfont\cyrfam=\twelvecyr
    \scriptscriptfont\cyrfam=\tencyr
  \setbox\strutbox=\hbox{\vrule
      height12.3pt depth5.54pt width0pt}%
  \baselineskip=16pt\rm}

\newbox\AuthorBox\newbox\TitleBox
\newbox\TFLinebox
\newbox\FLinebox
\newbox\HLinebox
\def\SetTFLinebox#1{\setbox\TFLinebox=\hbox{#1}}
\def\SetFLinebox#1{\setbox\FLinebox=\hbox{#1}}
\def\SetHLinebox#1{\setbox\HLinebox=\hbox{#1}}

 \def\SetAuthorHead#1{%
     \setbox\AuthorBox=\hbox{\ninepoint \it 
           \ignorespaces\frenchspacing#1\unskip}}
 \def\SetTitleHead#1{%
     \setbox\TitleBox=\hbox{\ninepoint \it
           \ignorespaces\frenchspacing#1\unskip}}

  \def\itSpacing{\relax}
  \def\itSpacingOff{\relax}


 \def\Hrule{\hrule width0pt height0pt}

  \newskip\ProcSkip \ProcSkip 8pt plus2pt minus2pt

 \newskip\LastSkip
 \def\SaveLastSkip{\LastSkip\lastskip}
 \def\RestoreLastSkip{\vskip-\LastSkip\vskip\LastSkip}

 \def\NoindentAfter{\everypar={\setbox0=\lastbox\everypar={}}}

 \long\def\H#1\par#2\par{\notenumber=0 \titlepagetrue%
    {
    \baselineskip=20pt
    \parindent=0pt\parskip=0pt\frenchspacing
    \leftskip=0pt plus .2\hsize minus .3\hsize
    \rightskip=0pt plus .2\hsize minus .3\hsize
 \def\\{\unskip\break}%
    \pretolerance=10000 \Hfont #1\unskip\break
     \vskip7pt\Hrule
\hfill \Authorfont #2\hfill\hfill\unskip}
    \vskip48pt plus 4pt minus 4pt
    \par\NoindentAfter\rm}

 \long\def\Hi#1\par#2\par{\notenumber=0 \titlepagetrue%
    {  \baselineskip=0pt  \parindent=0pt\parskip=0pt\frenchspacing
    \leftskip=0pt plus .2\hsize minus .3\hsize
    \rightskip=0pt plus .2\hsize minus .3\hsize
}
    \rm}


 \newdimen\PageRemainder
  \def\SetPageRemainder{
     \PageRemainder=\pagegoal
     \ifdim\PageRemainder=\maxdimen\PageRemainder=\vsize
     \else\advance\PageRemainder by -1\pagetotal\fi}

  \def\Rpt@{}\def\Rpt@@{}

  \long\def\HH#1\par{\par
  \SaveLastSkip\removelastskip\goodbreak
  \ifdim\LastSkip<30pt 
     \LastSkip 30pt
plus 3pt minus 2pt\fi
  \SetPageRemainder\advance\PageRemainder-\LastSkip
  \ifdim\PageRemainder<150pt
       \edef\Rpt@{remain = \the\PageRemainder\noexpand\\
                pagetotal=\the\pagetotal\noexpand\\
                           pagegoal=\the\pagegoal}%
          \fi
   \ifdim\PageRemainder<65pt 
       \ifdim\PageRemainder > 0pt
          \edef\Rpt@@{\noexpand\\
                      Had HH PageRemainder$<$\relax 65pt\noexpand\\
                      Hence forced break!}%
     \vskip 0pt plus .2\PageRemainder\eject 
    \fi\fi
    \vskip\LastSkip\Hrule 
    \pretolerance=10000\rightskip=0pt plus 3em
    \hangafter1 \hangindent=2.2em%
    \noindent
    \HHfont \unskip \Ednote{\Rpt@\Rpt@@}%
            \def\Rpt@{}\def\Rpt@@{}%
            \ignorespaces
            #1\par\rightskip=0pt\pretolerance=\StdPretolerance%
    \NoindentAfter
\tenpoint\rm%
     \medskip \vskip\ProcSkip}

  \long\def\HHH#1\par{\par%
  \SaveLastSkip\removelastskip\goodbreak
  \ifdim\LastSkip<\ProcSkip%
     \LastSkip\ProcSkip\fi
  \SetPageRemainder\advance\PageRemainder-\LastSkip
  \ifdim\PageRemainder<150pt
       \edef\Rpt@{remain = \the\PageRemainder\noexpand\\
                pagetotal=\the\pagetotal\noexpand\\
                           pagegoal=\the\pagegoal}%
       \fi
   \ifdim\PageRemainder<48pt  
        \ifdim\PageRemainder > 0pt
             \edef\Rpt@@{\noexpand\\
                      Had HHH PageRemainder$<$\relax48pt\noexpand\\
                      Hence forced break!}%
       \vskip 0pt plus .2\PageRemainder\eject 
      \fi\fi
   \vskip\LastSkip\par\noindent
   \HHHfont \unskip\Ednote{\Rpt@\Rpt@@}%
  \def\Rpt@{}\def\Rpt@@{}%
  \ignorespaces
   #1\unskip.\quad\rm\ignorespaces
   \ignorepars}

  \long\def\ignorepars#1\par{\def\Test{#1}%
     \ifx\Test\Empty\def\This{\ignorepars}%
        \else\def\This{\Test\par}\fi
           \This}
  \def\Empty{}

 \def\Abstract#1\par{\bgroup\Smallfonts\narrower\HHH #1\par}
 \def\endAbstract{\par\egroup}


 \def\ProcBreak{\par%
    \ifdim\lastskip<8pt%
    \removelastskip%
    \penalty-200\vskip\ProcSkip\fi}

 \def\th#1\par{\ProcBreak \noindent
   {\thfont\ignorespaces
    #1\unskip.}\it\itSpacing\kern.4em\ignorepars}

 \def\endth{\ProcBreak\rm\itSpacingOff }


 \def\pf#1\par{\ProcBreak %
    \noindent\pffont#1\unskip.\rm\itSpacingOff{\kern .7em}\ignorepars}

 \def\endpf{\medskip \ProcBreak } 

  \def\qedbox{\hbox{\vbox{
    \hrule width0.2cm height0.2pt
    \hbox to 0.2cm{\vrule height 0.2cm width 0.2pt
             \hfil\vrule height0.2cm width 0.2pt}
    \hrule width0.2cm height 0.2pt}\kern1pt}}

  \def\qed{\ifmmode\qedbox
    \else\unskip\ \hglue0mm\hfill\qedbox\ProcBreak\fi}

  \def \rk #1\par{\ProcBreak
     \noindent{\rkfont\ignorespaces #1\unskip.}%
     \rm\kern.6em\ignorepars}

  \def \endrk {\medskip\ProcBreak }

  \def \df #1\par{\ProcBreak
     \noindent{\dffont\unskip\ignorespaces #1\unskip.}%
     \rm\kern.6em\ignorepars}

  \def \enddf {\medskip\ProcBreak }

  \def \eg #1\par{\ProcBreak
     \noindent\egfont\unskip\ignorespaces #1\unskip.
     \rm\kern.6em\ignorepars}

  \newdimen\Overhang

   \def\MaxTag@#1#2#3#4#5{\setbox0=\hbox{#4\ignorespaces#2\unskip}%
     \dimen0=\wd0\advance\dimen0 by#3
     \ifdim\dimen0<#5\relax\dimen0=#5\fi
     \expandafter\edef\csname #1Hang\endcsname{\the\dimen0}}

 \def\MaxItemTag#1{\MaxTag@{Item}{#1}{.4em}{\ItemStyle}{\parindent}}%
 \def\MaxItemItemTag#1{%
        \MaxTag@{ItemItem}{#1}{.4em}{\ItemItemStyle}{\parindent}}
 \def\MaxNrTag#1{\MaxTag@{Nr}{#1}{.5em}{\NrStyle}{\parindent}}
 \def\MaxReferenceTag#1{%
        \MaxTag@{Reference}{[#1]}{.6em}{\ninerm}{\parindent}}
 \def\MaxFootTag#1{\MaxTag@{Foot}{#1}{.4em}{\ninerm}{\z@}}

  \def\SetOverhang@{\Overhang=.8\dimen0%
     \advance\Overhang by \wd0\relax
     \ifdim\Overhang>\hangindent\relax
       \advance\Overhang by .25\dimen0%
       \Ednote{Tag is pushing text.}\osumess{Tag is pushing text.}%
     \else\Overhang=\hangindent
     \fi}

   \def\Item#1{\par\noindent
      \hangafter1\hangindent=\ItemHang
      \setbox0=\hbox{\ItemStyle\ignorespaces#1\unskip}%
      \dimen0=.4em\SetOverhang@
      \rlap{\box0}\kern\Overhang\ignorespaces}

   \def\ItemItem#1{\par\noindent
      \hangafter1\hangindent=\ItemItemHang
      \setbox0=\hbox{\ItemItemStyle\ignorespaces#1\unskip}%
      \dimen0=.4em\SetOverhang@
      \advance\hangindent by \ItemHang
      \kern\ItemHang\rlap{\box0}%
      \kern\Overhang\ignorespaces}

  \def\Nr#1{\par\noindent\hangindent=\NrHang 
    \setbox0=\hbox{\NrStyle\ignorespaces#1\unskip}%
    \dimen0=.5em\SetOverhang@
    \rlap{\box0}\kern\Overhang
    \hangindent=\z@\ignorespaces}

   \newskip\Rosterskip\Rosterskip 1pt plus1pt 
   \def\Roster{\par\ifdim\lastskip<\Rosterskip\removelastskip\vskip\Rosterskip\fi
    \bgroup}
   \def\endRoster{\par\global\edef\LastSkip@{\the\lastskip}\removelastskip
       \egroup\penalty-50\LastSkip\LastSkip@\relax
       \ifdim\LastSkip<\Rosterskip\LastSkip\Rosterskip\fi
       \vskip\LastSkip}




 \def\cite#1{
    \def\nextiii@##1,##2\end@{{\frenchspacing\rm 
      \lBr\ignorespaces##1\unskip{\rm,~\ignorespaces##2}\rBr}}%
    \IN@0,@#1@%
    \ifIN@\def\next{\nextiii@#1\end@}\else
    \def\next{{\rm\lBr#1\rBr}}\fi\next}


   \def \Bib#1\par{%
       \par\removelastskip\SetPageRemainder
       \ifdim\PageRemainder < 97pt
        \ifdim\PageRemainder > 0pt
        \vfill\eject
       \fi\fi
    \ProcBreak \par\begingroup\parskip=0 pt%
    \goodbreak \vskip 15 pt plus 10 pt
    \noindent\null\hfill\Bibfont
      \ignorespaces #1\unskip\hfill\null\par 
    \frenchspacing \Smallfonts\rm
    \parskip=2.5 pt plus 1 pt minus.5pt%
    \nobreak\vskip 12pt plus 2pt minus2pt\nobreak
    \leftskip=0 pt \baselineskip=10.5pt}

 \def\ReferenceTagSlide{0em}
  \def\ReferenceTagGap{.5em}

  \def \rf#1{\par\noindent
     \hangafter1\hangindent=\ReferenceHang      
     \setbox0=\hbox{\ninerm[\ignorespaces#1\unskip]}%
     \dimen0=\ReferenceTagGap\SetOverhang@
     \rlap{\kern\ReferenceTagSlide\box0}%
     \kern\Overhang\ignorespaces}

  \def\ref#1\par#2\par#3\par#4\par{%
     \rf{#1}#2\unskip,\ #3\unskip,\
     #4\unskip.}

  \def\endBib{\par\endgroup\vskip 12pt minus 6pt }


  \long\def\Coordinates#1\endCoordinates{
 {\par\vskip4pt\def\\{\unskip, }\Coordfont\baselineskip10.5pt\noindent#1}}

 \def\pagecontents{
  \gdef\Pagetot@l{\pagetotal}
  \ifvoid\TRMargIns\else
    \rlap{\kern\hsize\kern10pt\vbox to 0pt{%
         \box\TRMargIns\vss}}\fi
  \ifvoid\topins\else\unvbox\topins\fi
   \dimen@=\dp\@cclv \unvbox\@cclv 
   \ifvoid\footins\else 
     \vskip\skip\footins
     \footnoterule
     \unvbox\footins\fi
   \ifr@ggedbottom \kern-\dimen@ \vfil \fi}


 \newcount\Ht 

 \def \Acc{\expandafter } 

 \def\swthat{\raise -1.1 ex\hbox{\sam$\widehat{}$}}
 \def\swttilde{\raise -1.2 ex\hbox{\sam$\widetilde{}$}}
 \def \overdot{{\raise .2 ex \hbox to 0pt {\hss\bf\smash{.}\hss}}}
 \def \overcircle{{\raise .1 ex \hbox to 0pt
    {\sam$\eightpoint\scriptstyle\hss\circ\hss$}}}

 \def \Mathaccent#1#2{{\sam 
  \setbox4=\hbox{$\vphantom{#2}$}
  \Ht=\ht4 
  \setbox5=\hbox{${#1}$}
  \setbox6=\hbox{${#2}$}
  \setbox7=\hbox to .5\wd6{}
  \copy7\kern .1\Ht \raise\Ht sp\hbox{\copy5}\kern-.1\Ht
  \copy7\llap{\box6}
  }}

  \def\SwtCheck #1{
        \ifmmode \check{#1}%
                \else \v {#1}%
                \fi}

 \def\barpartial {%
   \kern .17 em
    \overline {\kern -.17 em\partial\kern-.03 em}%
    \kern .03 em}

 
  \def\Overline#1{\setbox1=\hbox{\sam ${#1}$}%
      \ifdim \wd1 > 6pt
    \kern .11 em
    \overline {\kern -.11 em#1\kern-.14 em}
    \kern .14 em
  \else
    \kern .03 em
    \overline {\kern -.03 em#1\kern-.04 em}
    \kern .04 em
  \fi}

 \def\SOverline#1{\setbox1=\hbox{\sam ${#1}$}%
      \ifdim \wd1 > 7pt
    \kern .22 em
    \overline {\kern -.22 em#1\kern-.09 em}%
    \kern .09 em
  \else
    \kern .10 em
    \overline {\kern -.10 em#1\kern-.04 em}%
    \kern .04 em
  \fi}


 \def\Underline#1{\setbox1=\hbox{\sam ${#1}$}%
      \ifdim \wd1 > 6pt
    \kern .11 em
    \underline {\kern -.11 em#1\kern-.14 em}
    \kern .14 em
  \else
    \kern .03 em
    \underline {\kern -.03 em#1\kern-.04 em}
    \kern .04 em
  \fi}

 \def\SUnderline#1{\setbox1=\hbox{\sam ${#1}$}%
      \ifdim \wd1 > 7pt
    \kern .04 em
    \underline {\kern -.04 em#1\kern-.2 em}%
    \kern .2 em
  \else
    \kern .0 em
    \underline {\kern -.0 em#1\kern-.15 em}%
    \kern .15 em
  \fi}


 \def \Blackbox
   {\leavevmode\hskip .3pt \vbox
   {\hrule height 5pt\hbox{\hskip 4.5pt}}\hskip .5pt}

 \def \XX{\Blackbox\kern.5pt\Blackbox} 

  \def\.{.\kern1pt}

    \def\Hyphen{\edef\this{\the\hyphenchar\font}%
          \hyphenchar\font=-1\char\this\hyphenchar\font=\this}

 \ifx\undefined\text
  \def\text#1{\hbox{\rm #1}}\fi 



   \everymath{}  

  \def\PassMath@@{\aftergroup\AfterMath@} 

 \let\PassMath@\PassMath@@

 \def\AfterMath@{\futurelet\next\AfterMathMole@}

 \def\AfterMathMole@{
      \ifcat\next\space
          \def\this{}
      \else
      \ifcat\next\egroup %
        \def\this{\osumess{Handset mathsurround?? ---(see dollar brace)}}%
      \else
      \def\this{\AAfterMath@}
      \fi\fi
      \this}

 \def\hyphen@{-}
 \def\paren@{)}
 \def\apostr@{'}

 \def\MSC#1{\kern-.8\mathsurround#1\kern.8\mathsurround}

 \def\AAfterMath@#1{\def\Next{#1}
    \IN@0\Next @,.;:!?\relax @%
    \ifIN@\def\this{\MSC{\Next}}%
    \else
    \ifx\Next\hyphen@\def\this{\futurelet\next\AfterHyphen@}%
    \else
    \ifx\Next\paren@\def\this{#1}%
    \else 
    \ifx\Next\apostr@\def\this{#1}%
    \else \def\this{\osumess{Handset mathsurround??}%
                 #1}\fi\fi\fi\fi
    \this}

 \def\AfterHyphen@#1{\def\Next{#1}%
   \ifx\Next\hyphen@\def\this{--}\else
   \ifcat\next\space%
   \def\this{\kern-\mathsurround\kern.05em- \Next}\else
   \def\this{\kern-\mathsurround\kern.05em\Hyphen\Next}\fi\fi\this}

 \def\sam{\mathsurround=\z@\let\PassMath@\relax}  %
 \def\mas{\mathsurround=\StdMathsurround\let\PassMath@\PassMath@@}
 
 \def\Mas{\mathsurround=\StdMathsurround
                \everymath{\PassMath@}\let\PassMath@\PassMath@@}

 \def\m@th{\mathsurround=\z@\everymath{}}

 \def\m@@th{\mathsurround=\z@\everymath={}\let\m@th\relax}

\def\underbar#1{$\setbox\z@\hbox{#1}\dp\z@\z@
      \m@th \underline{\box\z@}$\relax}

\def\mathhexbox#1#2#3{\leavevmode
  \hbox{\m@@th$\m@th \mathchar"#1#2#3$}}

\def\dots{\relax\ifmmode\ldots\else$\m@th\ldots\,$\relax\fi}

\def\dotfill{\cleaders\hbox{\m@@th$\m@th \mkern1.5mu.\mkern1.5mu$}\hfill}
\def\rightarrowfill{$\m@th\mathord-\mkern-6mu%
  \cleaders\hbox{\m@@th$\mkern-2mu\mathord-\mkern-2mu$}\hfill
  \mkern-6mu\mathord\rightarrow$\relax}
\def\leftarrowfill{$\m@th\mathord\leftarrow\mkern-6mu%
  \cleaders\hbox{\m@@th$\mkern-2mu\mathord-\mkern-2mu$}\hfill
  \mkern-6mu\mathord-$\relax}

\def\downbracefill{$\m@th\braceld\leaders\vrule\hfill\braceru
  \bracelu\leaders\vrule\hfill\bracerd$\relax}
\def\upbracefill{$\m@th\bracelu\leaders\vrule\hfill\bracerd
  \braceld\leaders\vrule\hfill\braceru$\relax}

\def\angle{{\vbox{\m@@th\ialign{$\m@th\scriptstyle##$\crcr
      \not\mathrel{\mkern14mu}\crcr
      \noalign{\nointerlineskip}
      \mkern2.5mu\leaders\hrule height.34pt\hfill\mkern2.5mu\crcr}}}}

\def\big#1{{\m@@th\hbox{$\left#1\vbox to8.5\p@{}\right.\n@space$}}}
\def\Big#1{{\m@@th\hbox{$\left#1\vbox to11.5\p@{}\right.\n@space$}}}
\def\bigg#1{{\m@@th\hbox{$\left#1\vbox to14.5\p@{}\right.\n@space$}}}
\def\Bigg#1{{\m@@th\hbox{$\left#1\vbox to17.5\p@{}\right.\n@space$}}}
\def\n@space{\nulldelimiterspace\z@ \m@th}

\def\root#1\of{\setbox\rootbox\hbox{\m@@th$\m@th\scriptscriptstyle{#1}$}
  \mathpalette\r@@t}
\def\r@@t#1#2{\setbox\z@\hbox{\m@@th$\m@th#1\sqrt{#2}$\relax}
  \dimen@\ht\z@ \advance\dimen@-\dp\z@
  \mkern5mu\raise.6\dimen@\copy\rootbox \mkern-10mu \box\z@}

\def\mathph@nt#1#2{\setbox\z@\hbox{\m@@th$\m@th#1{#2}$}\finph@nt}

\def\mathsm@sh#1#2{\setbox\z@\hbox{\m@@th$\m@th#1{#2}$}\finsm@sh}

\def\@vereq#1#2{\lower.5\p@\vbox{\m@@th\baselineskip\z@skip\lineskip-.5\p@
    \ialign{$\m@th#1\hfil##\hfil$\crcr#2\crcr=\crcr}}}

\def\mathpalette#1#2{\sam\mathchoice{#1\displaystyle{#2}}%
  {#1\textstyle{#2}}{#1\scriptstyle{#2}}{#1\scriptscriptstyle{#2}}\mas}

\def\widehat#1{\setbox\z@\hbox{\sam$#1$}%
 \ifdim\wd\z@>\tw@ em\mathaccent"0\msbfam@5B{#1}%
 \else\mathaccent"0362{#1}\fi}
\def\widetilde#1{\setbox\z@\hbox{\sam$#1$}%
 \ifdim\wd\z@>\tw@ em\mathaccent"0\msbfam@5D{#1}%
 \else\mathaccent"0365{#1}\fi}

 \def\dots{\relax{}
  \ifmmode\def\thedots{\mdots@}\else\def\thedots{\tdots@}\fi %
  \thedots}

 \let\@ldeqno\eqno\let\@ldleqno\leqno
 \def\eqno{\everymath{}\@ldeqno} \def\leqno{\everymath{}\@ldleqno}

  \let\@ldeqalignno\eqalignno
  \def\eqalignno#1{\sam\@ldeqalignno{#1}\mas}
  \let\@ldeqalign\eqalign
  \def\eqalign#1{\sam\@ldeqalign{#1}\mas}

 \def\overrightarrow#1{\vbox{\m@th\ialign{##\crcr
      \rightarrowfill\crcr\noalign{\kern-\p@\nointerlineskip}
      $\hfil\displaystyle{#1}\hfil$\crcr}}}
 \def\overleftarrow#1{\vbox{\m@th\ialign{##\crcr
      \leftarrowfill\crcr\noalign{\kern-\p@\nointerlineskip}
      $\hfil\displaystyle{#1}\hfil$\crcr}}}
 \def\overbrace#1{\mathop{\vbox{\m@th\ialign{##\crcr\noalign{\kern3\p@}
      \downbracefill\crcr\noalign{\kern3\p@\nointerlineskip}
      $\hfil\displaystyle{#1}\hfil$\crcr}}}\limits}
 \def\underbrace#1{\mathop{\vtop{\m@th\ialign{##\crcr
      $\hfil\displaystyle{#1}\hfil$\crcr\noalign{\kern3\p@\nointerlineskip}
      \upbracefill\crcr\noalign{\kern3\p@}}}}\limits}

  \let\@ldmatrix\matrix
  \let\end@ldmatrix\endmatrix
  \def\matrix{\sam\@ldmatrix}
  \def\endmatrix{\end@ldmatrix\mas}
  \let\@ldgather\gather
  \let\end@ldgather\endgather
  \def\gather{\sam\@ldgather}
  \def\endgather{\end@ldgather\mas}
  \let\@ldalign\align
  \let\end@ldalign\endalign
  \def\align{\sam\@ldalign}
  \def\endalign{\end@ldalign\mas}
  \let\@ldaligned\aligned
  \let\end@ldaligned\endaligned
  \def\aligned{\sam\@ldaligned}
  \def\endaligned{\end@ldaligned\mas}
  \let\@ldtag\tag
  \def\tag{\sam\@ldtag}
   %

   \let\MinCDArrowWidth\minCDaw@




\newskip\insertskipamount\newskip\inserthardskipamount
\insertskipamount 6pt plus2pt 
\inserthardskipamount 6pt
\def\insertskip{\vskip\insertskipamount}
\newcount\SplitTest
\def\SetSplitTest{\SplitTest\insertpenalties
  \insert\topins{\floatingpenalty1}%
  \advance\SplitTest-\insertpenalties}
\def\midinsert{\par
 \SaveLastSkip\penalty-150\SetSplitTest\RestoreLastSkip
 \ifnum\SplitTest=-1
  \@midfalse\p@gefalse\else\@midtrue\fi\@ins}
\def\@ins{\par\begingroup\setbox\z@\vbox\bgroup%
  \vglue\inserthardskipamount}
\def\endinsert{\egroup 
  \if@mid \dimen@\ht\z@ \advance\dimen@\dp\z@
    \advance\dimen@\insertskipamount
    \advance\dimen@\pagetotal\advance\dimen@-\pageshrink
    \ifdim\dimen@>\pagegoal\@midfalse\p@gefalse\fi\fi
  \if@mid%
    \ifdim\lastskip<\insertskipamount\removelastskip\insertskip\fi
    \nointerlineskip\box\z@\penalty-200\insertskip
  \else%
    \SaveLastSkip
    \insert\topins{\penalty100 
    \splittopskip\z@skip
    \splitmaxdepth\maxdimen \floatingpenalty\z@
    \ifp@ge \dimen@\dp\z@
    \vbox to\vsize{\unvbox\z@\kern-\dimen@}
    \else \box\z@\nobreak\insertskip\fi}
    \RestoreLastSkip
   \fi\endgroup}


  \newcount\notenumber
  
  \def\note{\advance\notenumber by 1
    \footnote{\the\notenumber)}}

  \newbox\footbox

  \def\footnote#1{\let\@sf\empty
    \ifhmode\edef\@sf{\spacefactor\the\spacefactor}\/\fi
    \sam${}^{\fam0 #1}$\@sf\vfootnote{#1}}%

  \def\vfootnote#1{\insert\footins\bgroup
     \interlinepenalty100 \splittopskip=1pt
     \floatingpenalty=20000
     \leftskip=0pt\rightskip=0pt%
     \parindent=.3em
     \Smallfonts\rm
     \FootItem@{#1}
     \futurelet\next\fo@t}

  \def\FootItem@#1{\par\hangafter1\hangindent=\FootHang
     \setbox0=\hbox{\ignorespaces#1\unskip}%
     \dimen0=.4em\SetOverhang@
     \noindent\rlap{\box0}\kern\Overhang\ignorespaces}


  \def\fo@t{\ifcat\bgroup\noexpand\next \let\next\f@@t
    \else\let\next\f@t\fi \next}
  \def\f@@t{\bgroup\aftergroup\@foot\let\next}
  \def\f@t#1{\baselineskip=10pt\lineskip=1pt
            \lineskiplimit=0pt #1\@foot}%
  \def\@foot{
        \hbox{\vrule height0pt depth5pt width0pt}
        \egroup}
  \skip\footins=12 pt plus 0pt minus 0pt 
  \count\footins=1000 
  \dimen\footins=8in 



 \def\osumess#1{\EdSpider{\immediate\write16{Line \the\inputlineno: #1}}}%
 \def\HideEdStuff{\gdef\EdSpider##1{}}

 \font\BigSym=cmmi10 scaled \magstep 4

 \def\change{\InLMargin{\hbox{\BigSym \char63\kern10pt}}}

 \def\beginchange{\InLMargin{\hbox{\sam\twelvepoint$\heartsuit$\kern10pt}}}

 \def\endchange{\InLMargin{\hbox{\sam\twelvepoint$\spadesuit$\kern10pt}}}

 \def\InLMargin#1{\strut\vadjust{%
     \kern-\strutdepth
     \vtop to \strutdepth{%
         \baselineskip\strutdepth
         \llap{\sam$\smash{\hbox{\EdSpider{#1}}}$}\null}}}

 \def\strutdepth{\dp\strutbox}
 \def\strutheight{\ht\strutbox}

 \def\NoteInRMargin#1{\strut\vadjust{%
     \kern-1.001\strutdepth
     \vtop to \strutdepth{%
       \baselineskip\strutdepth
       \vss\rlap{\ninepoint\unskip\hskip\hsize
         \vtop to 0pt{%
           \hsize=16em\hfuzz=\hsize
           \leftskip=10pt%
           \rightskip=0pt plus 10000pt%
           \baselineskip=9.8pt\lineskip=.2pt%
           \let\\\break
           \noindent\EdSpider{#1}\vss}%
                \kern10pt}\hbox{}}
       }}

 \def\ednote#1{\NoteInRMargin{\tentt #1}}

 \def\cbar{\InLMargin{%
      \dimen0=\strutdepth\advance\dimen0 by \lineskip
      \vrule width 3pt
      height \strutheight depth \dimen0 \kern
      3pt}}

 \def\ccbar{\InLMargin{%
      \dimen0=2\strutdepth\advance\dimen0 by 2\lineskip
      \vrule width 3pt
        height 3\strutheight depth \dimen0 \kern
      3pt}}

 \newinsert\TRMargIns
 \dimen\TRMargIns=\maxdimen

  \def\Ednote#1{\insert\TRMargIns{%
       \vbox to 0pt{\hsize=140pt\hfuzz=\hsize
           \leftskip=6pt%
           \rightskip=0pt plus 10000pt%
           \baselineskip=9.8pt\lineskip=.2pt%
           \let\\\break
           \SetPageRemainder
           \vglue540pt\vglue-\PageRemainder
           \noindent\EdSpider{\tentt #1}\vss}%
       \smallskip}}

 \def\KillEdStuff{\def\ednote##1{}\def\Ednote##1{}%
      \let\change\relax\let\beginchange\relax\let\endchange\relax
       \let\cbar\relax\let\ccbar\relax}


  \topskip=12pt
  \newskip\StdBaselineskip 
  \StdBaselineskip 12pt
  \lineskip=1.1pt
  \lineskiplimit=.8pt
  \widowpenalty=10000 
  \clubpenalty=10000  
  \abovedisplayskip=6pt plus 1pt minus 1pt
  \abovedisplayshortskip=3pt plus 1.5pt
  \belowdisplayskip=6pt plus 1pt minus 1pt
  \belowdisplayshortskip=5pt plus 1pt minus 1pt
  \hfuzz=1.5pt   

  \def\StdPretolerance{100}
  \tolerance=\StdPretolerance

  \newdimen\StdMathsurround
  \StdMathsurround=1.5pt 
  \mathsurround=\StdMathsurround
  \Mas                   

   \def\prose{\relax\hbox{\kern.6\StdMathsurround}}
  
  \def\StdParskip{0pt}    
  \parskip=\StdParskip
  \parindent=0.5cm
 

  \def\Times{ptmr  } 
  \def\TimesI{ptmri  } 
  \def\TimesB{ptmb  }
  \def\TimesBI{ptmbi  }
  \def\HelveticaN{phvrrn }

  =\Times at 10bp
  =\TimesB at 10bp
  \font\tenit=\TimesI at 10bp
  =\TimesBI at 10bp

  \font\tenmrm=cmr10  


    =\Times at 9bp 
    \font\nineit=\TimesI at 9bp 
    =\TimesB at 9bp 
    =\TimesBI at 9bp 

    =\HelveticaN at 9bp 


  =\Times at 12bp
  \font\twelveit=\TimesI at 12bp
  =\TimesB at 12bp


  \font\titleit=\TimesI at 14.4bp
  =\TimesB at 14.4bp

 \SetAuthorHead{AuthorHead} 
 \SetTitleHead{TitleHead}  


  \def\lBr{\raise.125ex\hbox{[\kern.1125ex}}
  \def\rBr{\raise.125ex\hbox{\kern.1125ex]}}

 \setbox\footbox=\hbox{\Smallfonts 2)~}



  \bgroup
  \catcode`\@=11 
  \gdef\itSpacing{%
     \xspaceskip=.31em plus.1em minus.05em \sfcode `f=2001
     \itWarning@\let\itWarning@\itWarning@@}
  \gdef\itSpacingOff{%
     \xspaceskip=0pt \sfcode `f=1000
     \let\itWarning@\relax}
   \global\let\itWarning@\relax
  \gdef\itWarning@@{\errmessage{%
  Special italic spacing already in force
  (you have probably omitted an ``endth'').
  See itSpacing macro in osuPSfnt.sty
         }}
  \egroup

 \fontdimen1\titlebf=0.0pt
 \fontdimen2\titlebf=3.6135pt
 \fontdimen3\titlebf=2.8908pt
 \fontdimen4\titlebf=1.44539pt
 \fontdimen5\titlebf=6.64882pt
 \fontdimen6\titlebf=14.45398pt
 \fontdimen7\titlebf=1.60439pt

 \fontdimen1\tenbi=0.26794pt
 \fontdimen2\tenbi=2.50937pt
 \fontdimen3\tenbi=2.00749pt
 \fontdimen4\tenbi=1.00374pt
 \fontdimen5\tenbi=4.59717pt
 \fontdimen6\tenbi=10.03749pt
 \fontdimen7\tenbi=1.11415pt

 \fontdimen1\twelverm=0.0pt
 \fontdimen2\twelverm=3.01125pt
 \fontdimen3\twelverm=2.409pt
 \fontdimen4\twelverm=1.2045pt
 \fontdimen5\twelverm=5.39615pt
 \fontdimen6\twelverm=12.045pt
 \fontdimen7\twelverm=1.33699pt

 \fontdimen1\twelveit=0.27731pt
 \fontdimen2\twelveit=3.01125pt
 \fontdimen3\twelveit=2.409pt
 \fontdimen4\twelveit=1.2045pt
 \fontdimen5\twelveit=5.37207pt
 \fontdimen6\twelveit=12.045pt
 \fontdimen7\twelveit=1.33699pt

 \fontdimen1\twelvebf=0.0pt
 \fontdimen2\twelvebf=3.01125pt
 \fontdimen3\twelvebf=2.409pt
 \fontdimen4\twelvebf=1.2045pt
 \fontdimen5\twelvebf=5.5407pt
 \fontdimen6\twelvebf=12.045pt
 \fontdimen7\twelvebf=1.33699pt

 \fontdimen1\tenrm=0.0pt
 \fontdimen2\tenrm=2.50937pt
 \fontdimen3\tenrm=2.00749pt
 \fontdimen4\tenrm=1.00374pt
 \fontdimen5\tenrm=4.49678pt
 \fontdimen6\tenrm=10.03749pt
 \fontdimen7\tenrm=1.11415pt

 \fontdimen1\tenit=0.27731pt
 \fontdimen2\tenit=2.50937pt
 \fontdimen3\tenit=2.00749pt
 \fontdimen4\tenit=1.00374pt
 \fontdimen5\tenit=4.47672pt
 \fontdimen6\tenit=10.03749pt
 \fontdimen7\tenit=1.11415pt

 \fontdimen1\tenbf=0.0pt
 \fontdimen2\tenbf=2.50937pt
 \fontdimen3\tenbf=2.00749pt
 \fontdimen4\tenbf=1.00374pt
 \fontdimen5\tenbf=4.61723pt
 \fontdimen6\tenbf=10.03749pt
 \fontdimen7\tenbf=1.11415pt

 \fontdimen1\ninerm=0.0pt
 \fontdimen2\ninerm=2.25842pt
 \fontdimen3\ninerm=1.80673pt
 \fontdimen4\ninerm=0.90337pt
 \fontdimen5\ninerm=4.0471pt
 \fontdimen6\ninerm=9.03374pt
 \fontdimen7\ninerm=1.00273pt

 \fontdimen1\nineit=0.27731pt
 \fontdimen2\nineit=2.25842pt
 \fontdimen3\nineit=1.80673pt
 \fontdimen4\nineit=0.90337pt
 \fontdimen5\nineit=4.02904pt
 \fontdimen6\nineit=9.03374pt
 \fontdimen7\nineit=1.00273pt

 \fontdimen1\ninebf=0.0pt
 \fontdimen2\ninebf=2.25842pt
 \fontdimen3\ninebf=1.80673pt
 \fontdimen4\ninebf=0.90337pt
 \fontdimen5\ninebf=4.15552pt
 \fontdimen6\ninebf=9.03374pt
 \fontdimen7\ninebf=1.00273pt


 \newcount\MaxSpaceFactor
 \MaxSpaceFactor=3000 

 \def\ItemStyle{\rm}
 \def\NrStyle{\rm}
 \def\ItemItemStyle{\rm}

 \MaxItemTag{(iii)}
 \MaxItemItemTag{(iii)}
 \MaxNrTag{(2)}
 \MaxFootTag{2)}
 \def\ReferenceHang{30pt}

 \catcode`\@=\active


\loadbold

=\Times  
=\Times scaled750
=\Times scaled650
\font\rms=\Times scaled 920 

=\TimesBI scaled 860
=\TimesI scaled 860

\textfont0=\rrm  
\scriptfont0=\erm 
\scriptscriptfont0=\srm

\def\Augment#1#2{%
    \toks0\expandafter{#1}\toks2{#2}%
    \edef#1{\the\toks0\the\toks2}}

 \font\twelverma=\Times  scaled 1200
 \font\tenrma=\Times  scaled 1000
 \font\ninerma=\Times scaled 920
 =\Times scaled 840
 \font\sevenrma=\Times scaled 760
 =\Times scaled 680
 \font\fiverma=\Times scaled 600

 \Augment\tenpoint{%
  \textfont0=\tenrma  \scriptfont0=\sevenrma  
  \scriptscriptfont0=\fiverma  }

 \Augment\ninepoint{%
  \textfont0=\ninerma  \scriptfont0=\sevenrma 
  \scriptscriptfont0=\fiverma}

 \Augment\twelvepoint{%
  \textfont0=\twelverma  \scriptfont0=\ninerma  
  \scriptscriptfont0=\sevenrma}

\mathsurround=1pt
\hsize=13.45truecm
\vsize=19.5truecm
\hoffset=1.25truecm
\voffset=2truecm
\advance\baselineskip by 2pt

\predefine\til{\~}
\def\~#1{\relax\ifmmode\widetilde{#1}\else\til{#1}\fi}

\redefine \le{\leqslant}

\define \wt#1{\mathaccent"0365{#1}}
\define \wh#1{\mathaccent"0362{#1}}
\def\sdp{{\rtimes}}

\define \iss{\,\Mathaccent{\raise -.8 ex\hbox{$\widetilde{}$\kern.1em}}\rightarrow\,}

\define \bigcupp{{\overset\cdot\to\bigcup}}

\define \chr{\mathop{\fam0 char}\,}

\define \ind{\operatorname{\fam0 ind}}
\define \pro{\operatorname{\fam0 pro}}

\define \Hom{\operatorname{\fam0 Hom}}

\define \Stab{\operatorname{\fam0 Stab}}
\define \aff{\mathop{\fam0 aff}}

\Mas
\HideEdStuff
\rm 
 

\def\issn{{\nineit ISSN 1464-8997 (on line) 1464-8989 (printed)}}

\def\gtp{{\nineit Published 10 December 2000: \ \copyright\ Geometry \& 
Topology Publications}}

\def\gtv3{{\nineit Geometry \& Topology Monographs, Volume 3 (2000) --
Invitation to higher local fields}}


\def\lione
{{\rms Geometry \& Topology Monographs}}

\def \litwo{{\rms Volume 3: Invitation to higher local fields
}} 

\def\tinfo #1.#2.#3-#4
{{
\noindent  {\lione} \hfill 
\par 
\vskip-1.5pt
\noindent {\litwo} \hfill
\par 
\vskip-1,5pt
\noindent {\rms Part #1, section #2, pages #3--#4} \hfill
\vskip24pt 
}}

\def\tinfos #1.#2.#3-#4
{{
\noindent  {\lione} \hfill 
\par 
\vskip-1.5pt
\noindent {\litwo} \hfill
\par 
\vskip-1.5pt
\noindent {\rms Pages #3--#4} \hfill
\vskip24pt 
}}

\def\tinfoi #1
{{
\noindent  {\lione} \hfill 
\par 
\vskip-1.5pt
\noindent {\litwo} \hfill
\par 
\vskip-1.5pt
\noindent {\rms Pages iii--xi: Introduction and contents} \hfill
\vskip26pt 
}}


  \def\titlepagehead{\hfil}

  \newif\iftitlepage\titlepagefalse
  \newif\ifblankpage\blankpagefalse
  \def\makeheadline{
     \ifblankpage{}\else%
     \iftitlepage
\vbox{\line{\vbox to 8.5pt{}
\ninerm
\copy\HLinebox \hfill
\hglue5mm\ninebf\folio 
\titlepagehead}}%
      \else
\vbox{\ifodd\pageno\rightheadline\else\leftheadline\fi}%
      \fi\vskip 12pt\fi}%
     \def\rightheadline{\line{\vbox to 8.5pt{}%
      \ninerm
\copy\TitleBox \hfill
\hglue5mm\ninebf\folio}}%
     \def\leftheadline{\line{\vbox to 8.5pt{}%
        \unskip\ninerm\unskip\ninebf\folio\hglue5mm
 \hfill \copy\AuthorBox
}}

 \footline={\ifblankpage{}\else
\iftitlepage\ninepoint\sam\hfill
\line{\vbox to 8.5pt{}
\copy\TFLinebox
\hfill
\hglue5mm 
}
            \else
\ninepoint\sam\hfill
\line{\vbox to 8.5pt{}
\copy\FLinebox
\hfill 
\hglue5mm
}
\hfil\fi\global\titlepagefalse\fi}

\def\blankpage{{\blankpagetrue\noindent\hbox to 10pt{\hss}\vfill
\pagebreak}}

\tenpoint\rm 
 